\documentclass{amsart}

\newtheorem{theorem}{Theorem}[section]

\newtheorem{lemma}[theorem]{Lemma}

\begin{document}
\title{McKay's correspondence and characters of finite subgroups of SU(2)}
\author{W. Rossmann}
\address{Department \ of Mathematics,
University of Ottawa\\
Ottawa, Cananda}
\email{rossmann@uottawa.ca}
\maketitle
\begin{center}
Abstract
\end{center}

According to McKay \lbrack 1980\rbrack\ the irreducible characters of finite
subgroups of SU(2) are in a natural 1-1 correspondence with the extended
Coxeter-Dynkin graphs of type ADE. \ We show that the character values
themselves can be given by an uniform formula, as special values of
polynomials which arise naturally as numerators of Poincar\a'e series
associated to finite subgroups of SU(2) acting on polynomials in two
variables. These polynomials have been the subject of a number of
investigations, but their interpretation as characters has apparently not
been noticed.

\section{Introduction}

\noindent In 1980 McKay announced his astounding discovery that the finite
subgroups of \text{SU}$(2)$ are in natural 1-1 correspondence with the
extended Coxeter-Dynkin graphs of type ADE in the following way. Let $K$ be a
finite subgroup of \text{SU}$(2)$, $\{\chi_i\}$ its irreducible characters,
and $\chi$ the character of its natural representation on $\mathbb{C}^2$. Let
$M=(m_{ij})$ be the matrix defined by
\begin{equation*}
\chi\chi_i=\sum_jm_{ij}\chi_j.
\end{equation*}
\emph{The matrices }$M$\emph{ corresponding to the finite subgroups of}
\text{\emph{SU}}$(2)$\emph{ exactly the matrices of the form }$M=2I-C$\emph{
where }$C$\emph{ is the Cartan matrix of an extended Coxeter-Dynkin graph of
type ADE}. McKay apparently found and verified this fact by direct
computation. In the meantime there have been many attempts to explain it in
other ways or to provide further insight into this phenomenon. Steinberg
offered an explanation in terms of representation theory of finite groups in
1982. Gonzales-Sprinberg and J.-L. Verdier \lbrack 1983\rbrack\ gave an
explanation in terms of algebraic geometry, an approach also taken up by
Kn\"orrer \lbrack 1985\rbrack. \ Kostant \lbrack 1985\rbrack\ found a
remarkable relation between certain Poincar\a'e polynomials associated to
these groups and the action of the Coxeter element on the ADE root system,
and Springer \lbrack 1987\rbrack\ gave another method for the computation of
these polynomials. They play an important role in all papers mentioned,
except for Steinberg's. Mysteriously, they also appeared in an entirely
different context in a paper of Lusztig in \lbrack 1983\rbrack\ and
reappeared in his paper in \lbrack 1999\rbrack. It turns out that they have
another striking feature, also rather strange at first sight: these
polynomials 'are' the irreducible characters of the finite subgroups of
\text{SU}$(2)$, if this statement is taken with a grain of salt. The precise
formulation is given in Theorem 3.1, the proof in \S 4. Some of the
constructions in the proof have parallels in Galois theory of algebraic
number fields having as common origin a method of localization for Galois
groups at prime ideals (Lemmas 4.1,4.2): here the Galois group is $K$ and its
localization is a maximal abelian subgroup $T$. The ADE graphs of \S 2 encode
those properties the configuration of the these $T$s inside of $K$ which are
needed in Theorem 3.1 to identify the polynomials as characters. 

\bigskip

\section{The graphs}

\noindent Let $V$ be a two dimensional complex vector with a unitary inner
product. Let $G$ be a finite subgroup of its unitary group \text{U}$(V)$, and
let $Z=G\cap{\text{U}}(\mathbb{C})$ its intersection with the center of
\text{U}$(V)$. We shall be concerned with maximal abelian subgroups of $G$,
typically denoted $H$ with normalizer $N(H)$.

\begin{lemma}(a)Any two maximal abelian subgroups $H$ of $G$ intersect in $Z$
only.

(b)The groups $W=N(H)/H$ have order 1 or 2.\end{lemma}

\noindent This is immediate from the observation that any non-scalar element
of \text{U}$(V)$ has two distinct eigenvalues, but it will be useful to keep
in mind how this happens. Each maximal abelian subgroup $H$ of $G$
corresponds to a pair of lines $P,P^{-}$, its eigenspaces. The $P$s which
arise in this way are the singular lines, invariant under some non-scalar
element in $G$. Let $\mathcal{P}$ denote the set of these $P$s, a subset of
$\mathbb{P}(V)$. Two $H$s are $G$-conjugate if and only if the corresponding
(unordered) pairs $\{P,P^{-}\}$ are. A conjugacy class of $H$s corresponds to
either 1 or 2 orbits of $P$s, depending on whether $P$ and $P^{-}$ belong to
the same orbit or not. In the first case $N(H)/H$ has order 2 and
interchanges $P$ and $P^{-}$; in the second case $N(H)/H$ has order 1. The
map $G/H\,\rightarrow\,G/N(H)$ is 2$:1$ in the first case, 1:1 in the second.
In the first case, when the orbits of $P$ and $P^{-}$ coincide, we call this
orbit \emph{doubled}, in the second case, when they are different, we call
these orbits \emph{coupled.} The same terminology applies to the
corresponding class of $H$s. The decomposition of $G$ into conjugacy classes
then gives 
\begin{equation*}
|\bar{G}|=1+\frac{|\bar{G}|}{|N(\bar{H})|}(|\bar{H}\mathopen|-1)+\cdots
\end{equation*}
where $\bar{G}=G/Z,\bar{H}=H/Z,$ and the sum is extended over a complete set
of representatives for the conjugacy classes of $H$s. The equation is the
class equation for these groups. Its solution for the possible values of
$|\bar{H}|,|N(\bar{H})|,|\bar{G}|$ is goes back to the beginnings of the
theory, as does the determination of the groups themselves. 

The finite subgroups $\bar{G}$ of
\text{U}$(V)/{\text{U}}(\mathbb{C})\,\approx\,{\text{PSU}}(2)\,\approx\,{%
\text{SO}}(3)$ are the \emph{polyhedral groups} . We list them under three
headings ADE in the notation of \lbrack Coxeter, 1974, p15\rbrack.

\smallskip

\noindent\textbf{Type A}. $\bar{G}=(p)\approx\mathfrak{C}_p$ , cyclic of
order $p,p=1,2\cdots$\smallskip

\noindent\textbf{Type D}. $\bar{G}$ $=(p,2)\approx\mathfrak{D}_p$, dihedral
of order $2p.,$ $p=1,2\cdots$

\noindent\textbf{Type E}. 

Tetrahedral: $\bar{G}=(3,3,3)\approx\mathfrak{A}_4$, alternating of order
$12$,

Octahedral: $\bar{G}=(4,3,2)\approx\mathfrak{S}_4$, symmetric of order $24$,

Icosahedral: $\bar{G}=(5,3,2)\approx\mathfrak{A}_5$, alternating of order
$60$.

\smallskip

\noindent The notation $\bar{G}=(5,3,2)$, for example, means that there are
three singular orbits of this group on
$\mathbb{P}(V)\,\approx\,\mathbb{P}^1\,\approx\,\mathbb{S}^2$ with
stabilizers $\bar{H}$ of orders $5,3,2$, the three cyclic groups of rotations
about three axes of symmetries of the icosahedron passing through a vertex,
face, or edge.

Of particular interest are the finite subgroups of \text{SU}$(V)$, which are
known as binary polyhedral groups or Klein groups. We use the letters $K,T$
instead of $G,H$ for these. Except for odd cyclic groups, they are all
inverse images of polyhedral groups under the map
\text{SU}$(V)\rightarrow{\text{SO}}(3)$ with kernel $\{\pm 1\}$ which here
replaces the map \ \text{U}$(V)\rightarrow{\text{SO}}(3)$ with kernel
\text{U}$(\mathbb{C})$. 

It is convenient to use the isomorphism of \text{SU}$(2)$ with the group
U$(\mathbb{H})$ of quaternions of norm one to represent $K$ as a subgroup of
\text{U}$(\mathbb{H})$. For this purpose $\mathbb{H}$ is considered as a left
vector space of $\mathbb{C}$ via $\sqrt{-1}x=ix$ with $i^2=-1$ in
$\mathbb{H}$. $V$ is taken to be $\mathbb{H}$ equipped with this complex
structure. The action of \text{U}$(\mathbb{H})$ on $\mathbb{H}$ by right
multiplication gives its identification with \text{U}$(V)$, written $u\cdot
x=xu^{-1}$ in order to have \text{U}$(V)$ act on the left, as usual. Three
generators $e_A,e_B,e_C$ for $K$ may be chosen in the form
\begin{equation*}
e_A=e^{\frac{\pi}{p_A}J_A},\quad e_B=e^{\frac{\pi}{p_B}J_B},\quad
e_C=e^{\frac{\pi}{p_C}J_C}
\end{equation*}
for suitable quaternions $J_A,J_B,J_C$ satisfying $J^2=-1$ \lbrack Coxeter
1974, p68\rbrack. The three cyclic subgroups $T_A,T_B,T_C$ generated by
$e_A,e_B,e_C$ are the isotropy groups at three base points $P_A,P_B,P_C$ for
the three singular orbits $\mathcal{P}_A,\mathcal{P}_B,\mathcal{P}_C$ of $K$
on $\mathbb{P}(V)$. The graphs below are designed to encode some properties
the configuration of the singular lines $P_A,P_B,P_C$ in $V$ through the
configuration of their stabilizers $T_A,T_B,T_C$ inside of $K$, as explained
in the legend. 

\medskip

\begin{tabular}[t]{|c|c|}
\hline
\begin{tabular}[t]{cc}
{\small A}$_{\boldsymbol{2k}}$&{\small odd cyclic}\\
$\langle 2k+1\rangle$&$\mathfrak{C}_{2k+1}$
\end{tabular}%
&%
\begin{tabular}[t]{cc}
{\small A}$_{\boldsymbol{2k-1}}$&{\small even cyclic}\\
$\langle 2k\rangle=\langle
k,k,1\rangle$&$\mathfrak{C}_{2k}=\mathfrak{C}_k^\ast$
\end{tabular}%
\\
\hline
$
\begin{matrix}
+1&+2&\,\cdots\,&+k\\
0&&&|\\
-1&-2&\,\cdots\,&-k
\end{matrix}
$&$
\begin{matrix}
+1&+2&\,\cdots\,&+(k-1)\\
0&&&\ast\\
-1&-2&\,\cdots\,&-(k-1)
\end{matrix}
$\\
\hline
\end{tabular}%

\begin{tabular}[t]{|c|c|}
\hline
\begin{tabular}[t]{cc}
{\small D}$_{\boldsymbol{n}}\boldsymbol{,n=}${\small odd}&{\small odd
dihedral}\\
$\langle n,2,2\rangle$&$\mathfrak{D}_n^\ast$
\end{tabular}%
&%
\begin{tabular}[t]{cc}
{\small D}$_{\boldsymbol{n}}\boldsymbol{,n=}${\small even}&{\small even
dihedral}\\
$\langle n,2,2\rangle$&$\mathfrak{D}_n^\ast$
\end{tabular}%
\\
\hline
$
\begin{matrix}
&&&&\pm 1&\\
1&2&\,\cdots\,&n-3&\ast&\pm 1\\
&0&&&&
\end{matrix}
$&$
\begin{matrix}
&&&&1&\\
1&2&\,\cdots\,&n-3&\ast&1\\
&0&&&&
\end{matrix}
$\\
\hline
\end{tabular}%

\begin{tabular}[t]{|c|c|}
\hline
\begin{tabular}[t]{cc}
{\small E}$_{\boldsymbol{6}}$&{\small tetrahedral}\\
$\langle 3,3,2\rangle$&$\mathfrak{A}_4^\ast$
\end{tabular}%
&%
\begin{tabular}[t]{cc}
{\small E}$_{\boldsymbol{7}}$&{\small octahedral}\\
$\langle 4,3,2\rangle$&$\mathfrak{S}_4^\ast$
\end{tabular}%
\\
\hline
$
\begin{matrix}
&&1&&\\
\pm 1&\pm 2&\ast&\pm 2&\pm 1\\
0&&&&
\end{matrix}
$&$
\begin{matrix}
&&&1&&\\
1&2&3&\ast&2&1\\
0&&&&&
\end{matrix}
$\\
\hline
\end{tabular}%

\begin{tabular}[t]{|c|}
\hline
\begin{tabular}[t]{cc}
{\small E}$_{\boldsymbol{8}}$&{\small icosahedral}\\
$\langle 5,3,2\rangle$&$\mathfrak{A}_5^\ast$
\end{tabular}%
\\
\hline
$
\begin{matrix}
&&&&1&&\\
1&2&3&4&\ast&2&1\\
0&&&&&&
\end{matrix}
$\\
\hline
\end{tabular}%

\medskip

\noindent\textbf{Legend}. The first line lists the polyhedral type, the
second a symbol for $K$ like $\langle 5,3,2\rangle$ and a name for it as
abstract group like $\mathfrak{A}_5^\ast$; the star indicates the extension
by $\{\pm 1\}$ via \text{SU}$(2)\rightarrow{\text{SO}}(3)$. The cyclic case
of type A, will be omitted throughout, as it would only contribute some
awkward complications in terminology. Each graph consists of three
\emph{branches}, strings of \emph{nodes} joined at a \emph{central node}
labeled $\ast$. We use the letters $A,B,C$ as labels for the three branches
and various items attached to these.

The three branches of the graph correspond to the three singular orbits
$\mathcal{P}_A,\mathcal{P}_B,\mathcal{P}_C$ of $K$ on $\mathcal{P}$, its
three singular orbits on $\mathbb{P}(V)$. The isotropy groups are cyclic and
for an appropriate choice of base points $P_A,P_B,P_C$ they can be taken to
be the three cyclic subgroups $T_A,T_B,T_C$ corresponding to the three
generators $e_A,e_B,e_C$ of $K$ mentioned above. These three subgroups form a
complete system of representatives for the maximal abelian subgroups of $K$,
but redundant for odd D$_n$ and for E$_6$. 

The elements of $T_A$ $=$ $\{e_A^{\pm n}=\exp(\pm n\frac{\pi}{p_A}J_A),$ $n=$
$1,$ $\cdots,$ $p_A\}$ are represented by the nodes labeled $n$ or $\pm n,$
to be specified as $nA$ or $\pm nA$ to indicate the branch $A$, if necessary.
Each node, whether labeled $n$ or $\pm n$, represents two elements $\exp(\pm
n\frac{\pi}{p}J)$, except when these coincide. A branch on which the label is
$n$ the group $T$ represents a doubled class in the sense explained above;
branches on which the label is $\pm n$ come in pairs, say $A,B$, and the
corresponding groups $T_A,T_B$ are coupled, hence conjugate within $K$. The
two branches $A,B$ are then related by a symmetry of the graph, which occurs
only for odd D$_n$ and for E$_6$ and the four elements $\exp(\pm
n\frac{\pi}{p_A}J_{A\text{ }}),\exp(\pm n\frac{\pi}{p_B}J_B)$ lie in two
conjugacy classes, each containing two of the elements, namely the two
corresponding to the choices $(+,-)$ or $(-,+)$ of the ambiguous signs
$(\pm,\pm)$: opposite signs correspond to conjugate elements, equal signs to
non-conjugate elements. The inversion involution $[c]\mapsto[c^{-1}]$ on
conjugacy classes interchanges coupled classes and fixes all others, and is
therefore represented by the symmetry of the graph which interchanges coupled
nodes on symmetric branches. 

The central node labeled $\ast$ corresponds to the central element $-1$ of
$K$ inside of \text{SU}$(2)$ and lies in all $T$s. The node labeled $0$ may
stand for the identity element and with this interpretation it should be
thought of as attached to the ends of \emph{all} branches; but it has been
placed next to a particular node so as to produce the usual extended
Coxeter-Dynkin graph, if taken as attached to\emph{ that} node \lbrack
Bourbaki, 1968\rbrack. This extra node is added so that the graph may do
double duty as McKay's character graph, associated to the Cartan matrix in
the way explained in the introduction. In this interpretation the extra node
stands for the trivial character. Detailed verifications are omitted, except
for the following lemma, which justifies most of the rules given above.

\begin{lemma}Two coupled nodes $\pm n,\pm n$ represent two conjugacy classes.
Each pair with opposite signs represents the same class, each pair with equal
signs represents inverse classes.\end{lemma}

\noindent\begin{proof} The coupling on nodes can occur only for D$_n$ and for
E$_6$, where the assertion be verified directly. For E$_6$, when $K=\langle
3,3,2\rangle$, the three vertices can be taken to be \lbrack Coxeter, 1974,
p76\rbrack
\begin{equation*}
J_A=\frac{1}{\sqrt{3}}(i+j+k),\,J_B=\frac{1}{\sqrt{3}}(i-j+k),\,J_C=i.
\end{equation*}
One has $e_C=e^{\frac{\pi}{2}i}=i$ and $iJ_Ai^{-1}=-J_B$. Hence
$e_Ce_Ae_{_C}^{-1}=e_B^{-1}$. This shows that conjugation by $e_C$ implements
the coupling of the two branches $A,B$ as indicated in the graph. For D$_n$
the verification is similar.\end{proof}

\noindent\textbf{The reflection group }$\boldsymbol{K}'$. To start with, let
$K\subset G$ be any 'normal extension' of $K\subset{\text{SU}}(V)$ to a
finite group $G\subset{\text{U}}(V)$. Then $\bar{K}\subset\bar{G}$ is a
normal extension of polyhedral groups in \text{SO}$(3)$. (The few
possibilities are well-known \lbrack Coxeter, 1974, \S 7.1-7.3\rbrack, but
will not be needed here.) The three singular orbits of $K$ and $G$ on
$\mathbb{P}(V)$ depend only on $\bar{K}$ and $\bar{G}$. In particular, if
$\bar{G}=\bar{K}$ then $G$ has the same singular $P$s as $K$, hence the
maximal abelian groups $H$ of $G$ are extensions of the maximal abelian
groups $T$ on $K$. Generally, $G$ permutes the singular $P$s of $K$ on
$\mathbb{P}(V)$ and $G/K$ permutes the three $K$-orbits of these $P$s. The
$K$-orbits corresponding to branches of the same length are permuted among
themselves, so that this action of $G/K$ must be trivial except for D$_n$ and
E$_6$. On the other hand, if $G/K$ leaves a singular $K$-orbit invariant,
then this situation is represented in the form $sP=P$ for some $s\in G$ and
some $P$ in this orbit. Similar remarks apply to the action of $G/K$ on the
conjugacy classes of $T$s. Consider in particular the situation when one of
the three $T$s for $K$ is invariant by some $s\in G$, i.e. $sTs^{-1}=T$. In
that case $s$ permutes the two invariant lines $P,P^{-}$ of $T,$ leading to
two cases, as indicated.

\medskip

Case (1). $s\cdot P=P,s\cdot P^{-}=P^{-}.\,\,T:
\begin{pmatrix}
\lambda&0\\
0&\lambda^{-}
\end{pmatrix}
,\,$ $s=
\begin{pmatrix}
\mu&0\\
0&\mu^{-}
\end{pmatrix}
$.

Case (2) $s\cdot P=P^{-},s\cdot P^{-}=P$. \ $T:
\begin{pmatrix}
\lambda&0\\
0&\lambda^{-}
\end{pmatrix}
,\,$ $s=
\begin{pmatrix}
0&\sigma\\
1&0
\end{pmatrix}
$.

\medskip

\noindent Case (1) occurs whenever an element $s$ of $G$ leaves invariant a
singular line $P$ of $K$. If $s$ is furthermore a reflection, i.e. leaves
pointwise fixed a subspace of codimension 1, then exactly one of
$\mu,\mu^{-}$ must be $=1$. 

It is a remarkable fact that each finite subgroup $K$ of \text{SU}$(2)$ is
contained with index 2 in a subgroup of \text{U}$(2)$ generated by three
reflections of order 2 \lbrack Coxeter, 1974, p93\rbrack. The reflection
group associated to $K=\langle p_A,p_B,p_C\rangle$ in this way is denoted
$K'=\langle p_A,p_B,p_C\rangle'$. In this case the action of the two element
group $K'/K$ on the three singular orbits
$\mathcal{P}_A$,$\mathcal{P}_B$,$\mathcal{P}_C$ gives an involution on the
graph encoding these orbits. The following lemma identifies this action of
$K'/K$ as the inversion involution $[c]\mapsto[c^{-1}]$ on conjugacy classes.

\begin{lemma}$K'/K$ acts by inversion on the conjugacy classes in
$K$.\end{lemma}

\begin{proof} \ It suffices to show that $K'/K$ acts by inversion on the
classes of the three generators $e_A,e_B,e_C$ of the three cyclic groups. In
any case $K'/K$ must map the class of one of these generator to a class of
the same order. Thus only coupled end-nodes $\pm 1,\pm 1$ need be considered,
for odd D$_n$ and for E$_6$. Take the latter case, for example. 

From the construction in \lbrack Coxeter, 1974, \S 9.5\rbrack\ the group
$K'=\langle 3,3,2\rangle'$ has as the two element group $\langle
4,3,2\rangle/\langle 3,3,2\rangle$ as a quotient and the action of $K'/K$ on
the classes in $K=\langle 3,3,2\rangle$ agrees with that of $\langle
4,3,2\rangle/\langle 3,3,2\rangle$. Thus it remains to check that $\langle
4,3,2\rangle/\langle 3,3,2\rangle$ acts by inversion on the classes in
$\langle 3,3,2\rangle$. It follows from the construction of $\langle
4,3,2\rangle$ from $\langle 3,3,2\rangle$ in \lbrack Coxeter, 1974, \S
7.3\rbrack\ that $\langle 4,3,2\rangle$ has an element which interchanges the
two generators $e_A$ and $e_B$ of $\langle 3,3,2\rangle$. The classes
$e_A$and $e_B$ in $\langle 3,3,2\rangle$ are indeed inverses of each other,
as can be read off the E$_6$ graph and has been verified in the proof of the
preceding lemma.\end{proof}

\noindent\textbf{Four groups associated to a Schwarz triangl}e. We mention a
few facts about the groups under consideration, but the further development
is logically independent of these. To get an overview it seems best to start
with a \emph{Schwarz triangle} of the first kind, a triangle on a 2-sphere,
with vertices $A,B,C$ whose angles $\pi/p_A,\pi/p_B,\pi/p_C$ satisfy the
inequality
\begin{equation*}
\frac{1}{p_A}+\frac{1}{p_B}+\frac{1}{p_C}>1.
\end{equation*}
There are four groups associated to such a Schwarz triangle. They may be
described as follows. The groups in question are transformation groups on a
real or complex space $V_\mathbb{R}$ or $V_\mathbb{C}$, copies of
$\mathbb{R}^3$ or $\mathbb{C}^2$ . We denote them
$K_\mathbb{R},K'_\mathbb{R},K_\mathbb{C},K'_\mathbb{C}$. $K_\mathbb{R}$ and
$K_\mathbb{C}$ are finite subgroups of \text{O}$(V_\mathbb{R})$ and
\text{U}$(V_\mathbb{C})$, $K'_\mathbb{R}$ and $K'_\mathbb{C}$ are reflection
groups in \text{O}$(V_\mathbb{R})$ and \text{SU}$'(V_\mathbb{C})$ containing
$K_\mathbb{R}$ and $K_\mathbb{C}$ as subgroups of index two.
(\text{SU}$'(V_\mathbb{C})$ is the group generated by
\text{SU}$(V_\mathbb{C})$ together with $i$ and consists of unitary matrices
of determinant $\pm 1$.)

Realize the 2-sphere carrying the triangle inside the 3-space of quaternions
$x\in\mathbb{H}$ satisfying $\bar{x}=-x$. It is then given by the equations
$\bar{x}=-x,x\bar{x}=1$, which may be combined into $x^2=-1$. \ Write
$J_A,J_B,J_C$ for the vertices $A,B,C$ considered as quaternions.
$V_\mathbb{R}$ or $V_\mathbb{C}$ are realized in terms of quaternions as
indicated below, a quaternion $i$ satisfying $i^2=-1$ being required for
$V_\mathbb{C}$. \ The groups have generators $t=t_A,t_B,t_C$ and
$s=s_A,s_B,s_C$ which act on $x\in\mathbb{H}$ through certain elements
$e=e_A,e_B,e_C$ and $u=u_A,u_B,u_C$ as indicated. \ The generators are
determined by the triangle ($e_A=\exp(\pi
J_A/p_A),\,u_A=\exp(\frac{\pi}{2}J_{BC})=\,J_{BC}:=\,[J_B,J_C]/\|[J_B,J_C]%
\|$) and the relations satisfied by them can be extracted from \lbrack
Coxeter, 1974\rbrack. The data are summarized in the tables below. The labels
$(p_A,p_B,p_C)$ and $\langle p_A,p_B,p_C\rangle$ for these groups are those
of Coxeter. Cyclic permutations of $A,B,C$ in the relations written out are
understood.

\begin{tabular}[t]{|c|c|c|c|c|}
\hline
\multicolumn{3}{|c|}{$V_\mathbb{R}:\bar{x}=-x$}&generators&relations(%
\&cyclic)\\
\hline
{SO}$(V_\mathbb{R})$&$K_\mathbb{R}$&$(p_A,p_B,%
p_C)$&$t:exe^{-1}$&$t_A^{p_A}=(t_At_Bt_C)=1$\\
\hline
{O}$(V_\mathbb{R})$&$K'_\mathbb{R}$&$(p_A,p_B,p_C)'$&$s:u%
\bar{x}u^{-1}$&$s_A^2=1,(s_As_B$)$^{p_C}=1$\\
\hline
\end{tabular}%

\begin{tabular}[t]{|c|c|c|c|c|}
\hline
\multicolumn{3}{|c|}{$V_\mathbb{C}:\sqrt{-1}x:=ix$}&generators&relations
(\&cyclic)\\
\hline
{SU}$(V_\mathbb{C})$&$K_\mathbb{C}$&$\langle
p_A,p_B,p_C\rangle$&$t:xe^{-1}$&$t_A^{p_A}=(t_At_Bt_C)=-1$\\
\hline
{SU}$'(V_\mathbb{C})$&$K'_\mathbb{C}$&$\langle
p_A,p_B,p_C\rangle'$&$s:ixu^{-1}$&$s_A^2=1,(s_As_B)^{p_C}=-1$\\
\hline
\end{tabular}%

\section{Statement of the theorem}

\noindent We shall need some facts concerning complex reflection groups, to
be applied to the group $K'$. These are better explained in a more general
setting, as in \lbrack Bourbaki, 1968\rbrack. Let $V$ be any
finite-dimensional complex vector space, $G$ a finite subgroup of GL$(V)$
generated by reflections acting on $V$ on the left. Let $S$ be the algebra of
polynomial functions on $V$. The group $G$ acts on $S$ on the right via
$fg=f\circ g$. In addition there is the natural action of $\mathbb{C}^\times$
on $V$ and on $S$, giving a representation of $\mathbb{C}^\times\times G$ on
$S$. The $\mathbb{C}^\times$-stable subspaces $M$ are the graded ones and
have a formal $\mathbb{C}^\times$-character, the Poincar\a'e series
$P_M\in\mathbb{C}[\![t]\!]$ of $M$ \lbrack Bourbaki, 1968, p103\rbrack. If
the subspace $M$ is stable under $G$ as well, then it has a formal
$\mathbb{C}^\times\times G$ character
$P_M\in\mathbb{C}[\![t]\!]\otimes{\text{Ch}}(G)$: if the variable $t$ is
identified with the standard character of $\mathbb{C}^\times$ then
$P_M=\,\sum m_{ik}\chi_it^k$ where $m_{ik}$ is the multiplicity of the
irreducible character $\chi_i$ of $G$ in the space homogenous polynomials of
degree $k$ in $M$. $P_M$ can be viewed as a formal power series with
characters of $G$ as coefficients or as a sum of characters of $G$ with
formal power series as coefficients. The Poincar\a'e series $P_S$ of the
symmetric algebra S$(V^\ast)$ and $P_\Lambda$ of the exterior algebra
$\Lambda(V^\ast)$ satisfy $P_S(t)P_\Lambda(-t)=1$. Furthermore
$P_\Lambda(-t)=\sum(-1)^k\Lambda_V^kt^k$ is a polynomial and
$P_S(t)=1/P_\Lambda(-t)$ expanded as a formal power series in $t$. 

Let $R$ the subalgebra of $G$-invariants in $S$. It is a polynomial algebra
in $\dim(V)$ algebraically independent homogeneous generators, whose degrees
$d_i$ are uniquely determined up to order. Its Poincar\a'e series is
$P_R(t)=\prod(1-t^{d_i})^{-1}$, expanded as formal power series. For any
$x\in V$ let $\mathfrak{m}_x$ be the ideal of $R$ vanishing at $x$ or
equivalently on the orbit $Gx$. The particular ideal $\mathfrak{m}_0$ of $R$
vanishing at 0 is graded as is the ideal $S\mathfrak{m}_0$ it generates in
$S$. Let $F$ be a $\mathbb{C}^\times\times G$-stable complement of
$S\mathfrak{m}_0$ in $S$. Then $S\approx R\otimes F$, which implies that 
\begin{equation*}
S/S\mathfrak{m}_x\approx(R/R\mathfrak{m}_x)\otimes F\approx F
\end{equation*}
for any $x$ in $V$; it is \ always a $G$-isomorphism, but a
$\mathbb{C}^\times\times G$ isomorphism only for $x=0$. The formal
$\mathbb{C}^\times\times G$-character of $S=R\otimes F$ evidently satisfies
$P_S(t)=P_R(t)P_F(t)$. This equation can be written as
$P_S(t)^{-1}P_F(t)=P_R(t)^{-1}$ i.e. 
\begin{equation*}
(\sum(-1)^k\Lambda_V^kt^k)P_F(t)=\prod(1-t^{d_i}).
\end{equation*}
It is an equation for $P_F(t)$ in $\mathbb{C}[\![t]\!]\otimes{\text{Ch}}(G)$.
Specializing now to the case $G=K',V\approx\mathbb{C}^2$ this equation reads 
\begin{equation*}
(1-{\text{tr}}_Vt+{\det}_Vt^2)P_F(t)=(1-t^{d_1})(1-t^{d_2}).
\end{equation*}
\ The determinant $\det_V$ is $\pm 1$ on $K'$ and after restriction of
characters to $K$ it becomes $\det_V\equiv 1$. The multiplication by
\text{tr}$_V$ in \text{Ch}$(K)$ is given by McKay's matrix $M=(m_{ij})$ in
the basis $\{\chi_i\}$, i.e. the incidence matrix of the graph in \S 2. The
degrees $d_1,d_2$ may be found from the fact that $d_1d_2$ is the order of
$K'$ and $(d_1-1)+(d_2-1)$ the number of reflections \lbrack Bourbaki, 1968,
p110,111\rbrack\ or located in the tables of Shephard and Todd \lbrack
1954\rbrack. We list them here for reference.

\begin{tabular}[t]{|l|l|l|l|l|l|}
\hline
$K'$&A$_n$&D$_n$&E$_6$&E$_7$&E$_8$\\
\hline
$d_1$&$n+1$&$2n-2$&$8$&$12$&$20$\\
\hline
$d_2$&$2$&$4$&$6$&$8$&$12$\\
\hline
\end{tabular}%

\medskip

\noindent Expand $P_F(t)=\sum\chi_iP_i(t)$ in terms of the irreducible
characters $\chi_i$ of $K$. The equation for $P_F(t)$ becomes a system of
linear equations for the unknown polynomials $P_i(t)$ with coefficient matrix
$M(t):=[1-m_{ij}t+t^2]$:
\begin{equation*}
\sum_{j=0}^l(1-m_{ij}t+t^2)P_i(t)\chi_i=(1-t^{d_1})(1-t^{d_2})\chi_0.
\end{equation*}
It may be solved for the polynomials $P_i(t)$ in a mechanical fashion. (A few
lines of Maple code which will do the job can be found in \S 5).

The polynomials $P_i(t)$ have been tabulated in a number of places, for
example in the papers of Gonzales-Sprinberg and Verdier \lbrack 1983\rbrack\
and of Kostant \lbrack 1985\rbrack. A closed formula is given by Cramer's
rule:
\begin{equation*}
P_i(t)=\frac{M_i(t)}{\det M(t)}.
\end{equation*}
$M_i(t)$ is obtained from $M(t)$ by replacing the 'column' $\chi_j$ by the
'column' $(1-t^{d_1})(1-t^{d_2})\chi_0$. The determinant is $\det
M(t)=\det(1-(2I-C)t+t^2)$ which can be evaluated by the formula for the
characteristic polynomial of a Cartan matrix \lbrack Bourbaki, Ch.V, \S
6,\#3\rbrack. Another formula is given by Kostant \lbrack 1985, Theorem
1.11\rbrack, which is particularly interesting because of the way it brings
in the ADE root system and its Coxeter transformation: it gives the
polynomials in the form $P_i(t)=\sum_{\varphi\in\Phi_i}t^{n(\varphi)}$ with
$\varphi$ running over a set $\Phi_i$ of roots and equipped with a natural
length function $n(\varphi)$, both defined in terms of the Coxeter element in
the Weyl group. A different characterization of the$\,P_i(t)$ in terms of
Hecke algebras can be found in \lbrack Springer 1985, \S 8.2\rbrack; it has
its origin in \lbrack Lusztig, 1983\rbrack\ and plays a prominent role in
\lbrack Lusztig, 1999\rbrack. These polynomial polynomials are related to the
characters of $K$ as follows.

\begin{theorem}Let $\chi_i$ be the irreducible character associated to the
node $i$ on the ADE graph and let $e^{n\frac{\pi}{p}J}$ be the element of the
maximal abelian subgroup associated to the node $n$ on a branch. The value of
$\pi_i$ $:=\chi_i+\bar{\chi}_i$ at $e^{n\frac{\pi}{p}J}$ is
\begin{equation*}
\pi_i(e^{n\frac{\pi}{p}J})=P_i(e^{n\frac{\pi}{p}i}),
\end{equation*}
obtained by the substituting $t=e^{n\frac{\pi}{p}i}$ into the polynomial
$P_i(t)$\end{theorem}

\noindent We add some remarks. Since every element of $K$ is conjugate to an
element of some maximal abelian subgroup $T$ this formula gives all character
values. While $\pi_i\equiv\chi_i+\bar{\chi}_i$ itself is not irreducible, the
irreducible characters $\chi_i$ themselves can easily be extracted as well.
This is clear if $\bar{\chi}_i=\chi_i$, i.e. if the corresponding node on the
graph is fixed by the inversion involution, as must always be the case unless
the node belongs to a coupled pair labeled $\pm n,\pm n$ on graph. The cases
remaining occur only for odd D$_n\,$and for E$_6$. Then one can use the
following rule for the splitting of $P_i=P_i^{+}+P_i^{-}$ which induces the
splitting $\pi_i\equiv\chi_i+\bar{\chi}_i$ of the character. All polynomials
$P_i(t)$ are of the form
\begin{equation*}
P_i(t)=\sum_{-h_i<k<h_i}t^{h_i+k}
\end{equation*}
The splitting $P_i=P_i^{+}+P_i^{-}$ is given by splitting the sum
symmetrically at $k=0$.

\section{Proof of the theorem}

\noindent We start in a slightly more general setting. Let $V$ be a
2-dimensional unitary space, $G$ any finite subgroup of \text{U}$(V)$. Let
$S$ be the ring of polynomial functions on $V$, $L$ its quotient field of
rational functions. Let $K,R$ be the invariants of $G$ in $L,S$. Consider the
extension $S\supset R$ of integral domains. The classical ramification theory
of prime ideals in extensions of Dedekind domains \lbrack Hecke, 1923,
Kap.V\rbrack\ is not immediately applicable, since $L$ is not an algebraic
function field, having transcendence degree 2 over $\mathbb{C}$. But if we
pass from $V\approx\mathbb{C}^2$ to
$(V-\{0\})/\mathbb{C}^\times\approx\mathbb{P}^1$ by considering only
$\mathbb{C}^\times$-invariant (= homogeneous) ideals of $S$, then the
classical facts about the ramification of ideals in extensions of Dedekind
domains remain applicable as far as needed. $\,$The details are spelled out
in the following two lemmas.

\begin{lemma} (a)The field extension $L\supset K$ is Galois with Galois group
$G$. Let $\mathfrak{p}$ be a prime ideal in $R$. The prime ideals of $S$ over
$\mathfrak{p}$ are permuted transitively by $G$: they are of the form
$\mathfrak{P}_g=\mathfrak{P}g,\,g\in H\backslash G$ where $H=H(\mathfrak{P})$
is the subgroup of $G$ leaving $\mathfrak{P}$ invariant. 

(b)The field extension L$(\mathfrak{P})\supset K(\mathfrak{p})$ of quotient
fields for $S/\mathfrak{P}\supset R/\mathfrak{p}$ is Galois with Galois group
$H/I$, where $I=I(\mathfrak{P})\,$is the normal subgroup of $H$ fixing
$S/\mathfrak{P}$.\end{lemma}

\noindent\begin{proof}These assertions are general facts \lbrack Bourbaki,
1985, Ch.V, \S 2,\text{n}$^{\text{o}}2$; 1968, Ch.V, \S
5,\text{n}$^{\text{o}}5$.\rbrack.\end{proof}

\begin{lemma}Assume that $\mathfrak{p}$ is a non-zero prime ideal of $R$
which ramifies in $S$, i.e. $\mathfrak{p}=\mathfrak{P}\cap R$ for some prime
ideal $\mathfrak{P}$ of $S$ with $I(\mathfrak{P})\neq 1$. Then 

(a)$\mathfrak{P}=Sa,\,a\in V^\ast$ a non-zero linear form. 

(b)$H$is the subgroup of $G$ leaving the subspace $P=\{a=0\}$ invariant, $I$
its \ \ cyclic normal subgroup leaving $P$ pointwise fixed.

(c) $S\mathfrak{p}$ decomposes as a product of ideals in the form 
\begin{equation*}
S\mathfrak{p}=\prod_{g\in H\backslash G}\mathfrak{P}_g^e
\end{equation*}
where $e$ is the order of \ $I$ . 

(d)The residue ring $S/\mathfrak{p}$ $(:=S/S\mathfrak{p})$ decomposes as a
direct product of rings 
\begin{equation*}
S/\mathfrak{p}\approx\prod_{g\in H\backslash G}S/\mathfrak{P}_g^e,
\end{equation*}
and this decomposition exhibits the representation of
$\mathbb{C}^\times\times G$ on $S/\mathfrak{p}$ as induced by the
representation of $\mathbb{C}^\times\times H$ on $S/\mathfrak{P}^e$
.\end{lemma}

\noindent\begin{proof}The assumption $I(\mathfrak{P})\neq 1$ implies that the
zero set $P\,$of $\mathfrak{P}$ in $V$ is pointwise fixed by a linear
transformation $\neq 1$. Since $V$ is 2-dimensional it must be of the form
$P=\{a=0\}$ for some linear form $a\in V^\ast$ and hence $\mathfrak{P}=Sa$.
This proves (a), and (b) is then clear \lbrack Bourbaki,1968, Ch.V, \S
5,\text{n}$^{\text{o}}5$\rbrack.

The ideal $\mathfrak{P}=Sa$ is $\mathbb{C}^\times$-invariant and maximal
among such ideals since $\dim V=2$. In particular, any two distinct ideals
$\mathfrak{P}_g\,,\mathfrak{P}_{g'}$ satisfy
$\mathfrak{P}_g\,+\mathfrak{P}_{g'}\,=S$. This implies the decompositions (c)
and (d) of $S\mathfrak{p}$ and of $S/\mathfrak{p}$ \lbrack Bourbaki, 1961,
Ch.II, \S 1,\text{n}$^{\text{o}}2$\rbrack. The last assertion then follows
from the definition of an induced representation \lbrack Serre, 1977,
p28\rbrack.\end{proof}

\noindent As they stand, the representations in the lemma are infinite
dimensional: the representation spaces $S/\mathfrak{p}$ and
$S/\mathfrak{P}_g^{e_g}$ are $\,$in fact modules over the group algebra with
coefficients in the ring $R/\mathfrak{p}$. Finite-dimensional representations
can be obtained as follows. In the setting of the lemma, part (b), let
$\mathfrak{m}$ be a maximal ideal of $R$ containing $\mathfrak{p}$. It is the
ideal of $R$ vanishing on the orbit $Gx$ for some $x\,$in the 1-dimensional
subspace $P=\{a=0\}$. Assume $x\neq 0$. Let $\mathfrak{M}_g$ be the maximal
ideal of $S$ vanishing at $gx$. These ideals are no longer invariant under
$\mathbb{C}^\times$. Consider however the group of elements $h\in G$
satisfying $hx=\lambda(h)x$ for some $\lambda(h)\in\mathbb{C}^\times$. Since
$P$ is 1-dimensional, this group is the subgroup $H$ of $G$ leaving $P$
invariant. Let $I$ be the subgroup of $H$ leaving $P$ pointwise fixed. Then
$\lambda$ $:H/I\rightarrow\mathbb{C}^\times$ is a multiplicative character of
$H/I$. Its image is a cyclic subgroup $C$ of $\mathbb{C}^\times$.

The following corollary (which remains valid with $e=1$ if $\mathfrak{m}$ is
unramified, i.e. not fixed by any nontrivial element of $G$) will give the
desired reduction to finite dimensions.

\begin{lemma}(a)The residue ring $S/\mathfrak{m}$ decomposes as a direct
product 
\begin{equation*}
S/\mathfrak{m}=\prod_{g\in I\backslash G}S/\mathfrak{M}_g^e,\tag*{(1)}
\end{equation*}
and this decomposition exhibits the representation of $C\times G$ on
$S/\mathfrak{m}$ as induced by the representation of $C\times I$ on
$S/\mathfrak{M}^e$. 

(b)The representation of $I$ on $S/\mathfrak{M}^e$ is its regular
representation: 
\begin{equation*}
S/\mathfrak{M}^e\approx\prod_{s\in I}\mathbb{C}_s.\tag*{(2)}
\end{equation*}
(c)The representation of $G$ on $S/\mathfrak{m}$ is 
\begin{equation*}
S/\mathfrak{m}\approx\sum_{\nu=0}^{e-1}{\text{\textrm{Ind}}}_I^G\epsilon^\nu.%
\tag*{(3)}
\end{equation*}
where $\epsilon\,$ is the character of $I$ on the line $\mathbb{C}a$ in
$V^\ast$.\end{lemma}

\noindent\begin{proof}(a)Consider the representation $S\mathfrak{m}$ as an
intersection of primary ideals in $S$. These are necessarily powers of the
maximal ideals $\mathfrak{M}_g$ of $\,S$ vanishing at the points $xg,\,g\in
G$. Since $S\mathfrak{p}=\prod_{g\in H\backslash G}\mathfrak{P}_g^e$ and
since $I$ is the subgroup leaving $\mathfrak{M}$ invariant this decomposition
is given by $S\mathfrak{m}=\prod_{g\in I\backslash G}\mathfrak{M}_g^e$. The
direct product decomposition of $S/\mathfrak{m}$ then follows, as does its
identification as an induced representation.

(b) The group $I$ consists of reflections in the subspace $P=\{a=0\}$ and may
be considered as a reflection group on the quotient space $V/P$. The
assertion then follows from general facts about reflection groups \lbrack
Bourbaki 1968, p107\rbrack. More directly, the$\,$ classes of
$1,a,\cdots,a^{e-1}$ form a basis for $S/\mathfrak{M}^e$ and provide the
decomposition of $S/\mathfrak{M}^e$ as a sum over the $e$ irreducible
characters of the cyclic group $I$ characteristic of the regular
representation$:$
\begin{equation*}
S/\mathfrak{M}^e\approx\sum_{\nu=0}^{e-1}\mathbb{C}a^\nu\tag*{(4)}
\end{equation*}
The components in this sum are the powers $\epsilon^\nu,\,$
$\nu=\,0,\,1,\,\cdots\,,\,e-1,\,$of the character $\epsilon\,$ by which $I$
acts on the line $\mathbb{C}a$ in $V^\ast$. 

(c)Follows from (2) and (4). \end{proof}

\noindent Now suppose $G\subset{\text{U}}(V)$ is a normal extension of
$K\subset{\text{SU}}(V)$ as at the beginning of \S 2. Let $T$ be one of the
three subgroups $T_A,T_B,T_C$ of $K$ and take for $x$ an eigenvector of $T$.
The subgroup $H$ of $G$ satisfying $xh=\lambda(h)x$ is one of the three
maximal abelian subgroups of $G$, namely the one which contains $T$. If the
subgroup $I$ of $G$ fixing $x$ is non-trivial, then it is cyclic and its
generator is a reflection $s$ which centralizes $T$. The decomposition
$V=P^{+}+P^{-}$ of $V$ into eigenspaces of $T$ produces the situation of case
(1) in \S 2:
\begin{equation*}
\,T:
\begin{pmatrix}
\lambda&0\\
0&\lambda^{-}
\end{pmatrix}
,\quad s=
\begin{pmatrix}
1&0\\
0&\epsilon
\end{pmatrix}
.
\end{equation*}
Here $x$ belongs to $P^{+}$ and now $\epsilon=\det|I$. The decomposition (3)
becomes
\begin{equation*}
S/\mathfrak{m}\approx\sum_{\nu=0}^{e-1}{\det}^\nu\otimes{\text{Ind}}_I^G1,%
\quad\,{\text{Ind}}_I^G1\approx\prod_{g\in I\backslash
G}S/\mathfrak{M}_g\tag*{(5)}
\end{equation*}
the second isomorphism coming from $S/\mathfrak{M}_g\approx\mathbb{C}$. This
equation remains also valid when $I=1$. \ Furthermore, under the homomorphism
$\lambda:H\rightarrow C$ the action of $C$ by scalar multiplications agrees
with that of $H$ by permutations of cosets $Ig\mapsto hIg=Ihg$, but from the
opposite side of $G$. This follows from the relation $xhg=\lambda(h)xg$. 

Assume further that $G$ is generated by reflections. Then $S$ is free over
$R$, say $S\approx R\otimes F$ as a tensor product over $\mathbb{C}$. Then
$S/\mathfrak{m}\approx(R/\mathfrak{m})\otimes F$ and
$R/\mathfrak{m}\approx\mathbb{C}$ leads to
\begin{equation*}
S/\mathfrak{m}\approx F.
\end{equation*}
If the subspace $F$ of $S$ is chosen to be stable under
$\mathbb{C}^\times\times G$ then this is a $\mathbb{C}^\times\times
G$-isomorphism and the properties concerning the representation of $C\times
G$ on $S/\mathfrak{m}$ carry over to $F$. 

At this point we return to the group $K'$ in the role of $G$. Thus $S\approx
R\otimes F$ and the representation of $K'$ on $F$ is equivalent to its
regular representation in the direct product of copies $\mathbb{C}_g$ of
$\mathbb{C}$,
\begin{equation*}
F\approx\prod_{g\in K'}\mathbb{C}_g.\text{ }\tag*{(6)}
\end{equation*}
The elements of the direct product are functions $f_g$ of $g\in K'$ on which
$K'$ acts by translations in the variable $g$ from the left: $f_g\mapsto
f_{l^{-1}g}$ . This representation of $K'$ extends to the biregular
representation of $K'\times K'$ by left-right translations $f_g\mapsto
f_{l^{-1}gr}$.

\begin{lemma}The biregular representation of $K'\times K'$ decomposes on
$K\times K$ in the form
\begin{equation*}
\sum\chi_i\otimes\pi_i,\quad\pi_i:=\bar{\chi}_i+\chi_i.\tag*{(7)}
\end{equation*}
The sum runs over the irreducible characters $\chi_i$ of $K$.\end{lemma}

\noindent\begin{proof} Consider the decomposition of the biregular character
of $K'\times K'$ into two pieces corresponding to the decomposition of
functions on $K'=K\cup sK$ into functions supported on one of the two cosets.
That the functions supported on $K$ give the regular representation $\sum_i$
$\chi_i\otimes\bar{\chi}_i$ is clear. That those on $sK$ give
$\sum\chi_i\otimes\chi_i$ follows from the fact that $K'/K$ acts by
$\chi\mapsto\bar{\chi}$ on the characters of $K$: the map $f_{_k}\mapsto
f_{sk}$ sends function on $K$ to functions on $sK,$ transforming the left
action $f_{lk}$ of $l\in K$ into the twisted action $f_{s(s^{-1}\text{
}ls)k}$.\end{proof}

\noindent\textbf{Remark}. The decomposition of the lemma can also be
understood as follows. An irreducible character $\chi'$ of $K'$ either
remains irreducible on $K$ or else decomposes into two a pair of characters
conjugate under $K'/K$. All irreducible characters of $K$ occur twice in this
way: if $\chi'=\bar{\chi}'$ then twice in this $\chi'$; if
$\chi'\neq\bar{\chi}'$ then once in $\chi'$ and once in $\bar{\chi}'$. These
facts follow from general properties of subgroups of index 2 \lbrack Weyl,
1939, p.159\rbrack. They can also be found in \lbrack Gonzales-Sprinberg and
Verdier, 1983\rbrack.\medskip

\noindent Compare the representation (5) of $C\times K'$ on
$S/\mathfrak{m}\approx F$ with the decomposition (6) of the biregular
representation of $K'\times K'$. In (5) the group $K'$ acting on
$S/\mathfrak{m}\approx F$ corresponds to the subgroup $K'\times 1$ of
$K'\times K'$ acting on $\prod_{g\in K'}\mathbb{C}_g$. The character of $K$
in (5) corresponds to character of $K\times 1$ in (7), obtained by evaluating
the right factors $\pi_i$ in the lemma at the identity element. The action of
$1\times K$ in (7) is not directly visible in (5). However, the action of the
subgroup $1\times T$ \emph{is} visible in (5) and corresponds to the action
of the subgroup $C$ of $\mathbb{C}^\times$ via the homomorphism
$\lambda:T\rightarrow C$, as remarked after (5). The group $T$ consists of
elements $e^{n\frac{\pi}{p}J}$ where $J\in\mathbb{H}$ satisfies $J^2=-1$. Its
eigenspaces are those of $J$ acting on $V\approx\mathbb{H}$ by right
multiplication with eigenvalues $\pm i$. Choose $x$ in the $+i$ eigenspace,
so that $xJ=ix$. Then $xe^{n\frac{\pi}{p}J}=e^{n\frac{\pi}{p}i}x$, i.e.
$\lambda(e^{n\frac{\pi}{p}J})=e^{n\frac{\pi}{p}i}$. The result is that under
the isomorphism (6) the character of $K\times T$ coming from the biregular
representation of $K'\times K'$ on $\prod_{g\in K'}\mathbb{C}_g$ agrees with
the character of $K\times C$ coming from the character of
$K\times\mathbb{C}^\times$ on $F$, i.e. 
\begin{equation*}
\sum P_i(e^{n\frac{\pi}{p}i})\chi_i=\sum\pi_i(e^{n\frac{\pi}{p}J})\chi_i.
\end{equation*}
Comparison of the coefficients of $\chi_i$ concludes the proof of the theorem.

\section{An example}

\noindent We list the polynomials $P_i(t)$ and the character values
$\chi_i(e^{n\frac{\pi}{p}J})$ for the binary icosahedral group $K=\langle
5,3,2\rangle$ of type E$_8$. The labels $1A,2A,3A,4A,\ast$ refer to the 5
nodes on branch $A$, the node $\ast$ being common to all branches. 

\medskip

$P_{1A}(t)=t+t^{11}+t^{19}+t^{29}$

$P_{2A}(t)=t^2+t^{10}+t^{12}+t^{18}+t^{20}+t^{28}$

$P_{3A}(t)=t^3+t^9+t^{11}+t^{13}+t^{17}+t^{19}+t^{21}+t^{27}$

$P_{4A}(t)=t^4+t^8+t^{10}+t^{12}+t^{14}+t^{16}+t^{18}+t^{20}+t^{22}+t^{26}$

$P_\ast(t)=\sum_{k=0}^5t^{15-2k}+t^{15+2k}$

$P_{2B}(t)=t^6+t^8+t^{12}+t^{14}+t^{16}+t^{18}+t^{22}+t^{24}$

$P_{1B}(t)=t^7+t^{13}+t^{17}+t^{23}$

$P_{1C}(t)=t^6+t^{10}+t^{14}+t^{16}+t^{20}+t^{24}$

$P_0(t)=1+t^{30}$

\medskip

\noindent Maple will compute these polynomials if instructed this way:

\medskip

$>$with(linalg):

$>$C := matrix(\lbrack\lbrack
2,-1,0,0,0,0,0,0,0\rbrack,\lbrack-1,2,-1,0,0,0,0,0,0\rbrack,\lbrack
0,-1,2,-1,0,0,0,0,0\rbrack, 

\lbrack 0,0,-1,2,-1,0,0,0,0\rbrack, \lbrack
0,0,0,-1,2,-1,0,0,0\rbrack,\lbrack 0,0,0,0,-1,2,-1,0,-1\rbrack,

\lbrack 0,0,0,0,0,-1,2,-1,0\rbrack, \lbrack 0,0,0,0,0,0,-1,2,0\rbrack,\lbrack
0,0,0,0,0,-1,0,0,2\rbrack\rbrack);

$>$A := 1+t{\char94}2-t*(2-C);

$>$B := ((1-t{\char94}12)*(1-t{\char94}20))*vector(\lbrack
1,0,0,0,0,0,0,0,0\rbrack);

$>$P := linsolve(A,B);

$>$sort(expand(P\lbrack 4\rbrack))

\medskip

\noindent The characters are $\pi_i=2\chi_i$ and the values
$\chi_i(e^{n\frac{\pi}{}J})$ are listed in the following table. Both $i=0$
and $n=0$ are omitted. Notation: 

\medskip

$\tau:=\frac{1+\sqrt{5}}{2}=e^{\frac{\pi{\text{i}}}{5}}+e^{-\frac{\pi{%
\text{i}}}{5}}=1+e^{\frac{2\pi{{\text{i}}}}{5}}+e^{-\frac{2\pi{{%
\text{i}}}}{5}},$

$\tau^{-}:=\frac{1-\sqrt{5}}{2}=e^{\frac{3\pi{{\text{i}}}}{5}}+e^{-\frac{3%
\pi{{\text{i}}}}{5}}=1+e^{\frac{4\pi{\text{i}}}{5}}+e^{-\frac{4\pi{%
\text{i}}}{5}},$

$[(x-\tau)(x+\tau^{-})=x^2-x-1]$.

\medskip

\begin{tabular}[t]{|c|c|c|c|c|c|}
\hline
&$e^{\pm\frac{\pi}{5}J_A}$&$e^{\pm\frac{3\pi}{5}J_A}$&$e^\ast=-1$&$e^{n\frac{%
\pi}{3}J_B}$&$e^{\frac{\pi}{2}J_C}$\\
\hline
$\chi_{1A}$&$\tau$&$\tau^{-}$&$-2$&$(-1)^{n+1}$&$0$\\
\hline
$\chi_{2A}$&$\tau$&$\tau^{-}$&$+3$&$0$&$-1$\\
\hline
$\chi_{3A}$&$1$&$1$&$-4$&$\pm 1$&$0$\\
\hline
$\chi_{4A}$&$0$&$0$&$+5$&$-1$&$0$\\
\hline
$\chi_\ast$&$-1$&$0$&$+6$&$0$&$i$\\
\hline
$\chi_{1B}$&$\tau^{-}$&$\tau$&$-2$&$(-1)^{n+1}$&$0$\\
\hline
$\chi_{2B}$&$1$&$-1$&$+4$&$1$&$0$\\
\hline
$\chi_{1C}$&$\tau^{-}$&$\tau$&$+3$&$0$&$-1$\\
\hline
\end{tabular}%

\section{Remarks on invariants and geometry}

\noindent\textbf{Invariants}. The proof of the theorem was given in the
abstract setting of ramification of prime ideals, which seemed clearest and
perhaps useful in other situations. It is possible, however, to give a
description of the invariants of a finite subgroup $G$ of \text{U}$(2)$ which
makes the ring $R$ and the ramification of its homogenous prime ideals in $S$
more explicit. Let $P$ be the 1-dimensional subspace of
$V\approx\mathbb{C}^2$, and let $a$ be a linear form defining it as
$P=\{a=0\}$. Let $H$ and $I$ the subgroups of $G$ leaving $P$ invariant and
pointwise fixed. If $h\in H$ then $ah=\lambda_P(h)a$ for some multiplicative
character $\lambda_P:H\rightarrow\mathbb{C}^\times$. \ Let 
\begin{equation*}
f_P(x):=\prod_{g\in[H\backslash G]}a(gx),
\end{equation*}
the product being taken over a set of coset representative denoted
$[H\backslash G]$. Under the left action $x\mapsto lx$ by $l\in G$ the
polynomial $f_P(x)$ transforms according to the rule 
\begin{equation*}
f_P(lx)=\mu_P(l)f_P(x),\quad\mu_P:=(\lambda_P)^{|H\backslash G|},
\end{equation*}
showing that $f_P$ is a relative invariant of $G$. The zero-set $\{f_P=0\}$
is $GP$, the union of the lines $P_g:=gP,g\in G/H$. Fix an affine line
$x=az+b$ in $V$ and let $z_P$ be its point of intersection with $P$. The
$P_g$ then become the roots $z_g$ of the equation $f_P(z)=0$ which represents
$f_P=0$. Now consider any relatively invariant form $f$ on $V$. The equation
$f=0$ defines a set of points on $\mathbb{P}(V)$ invariant under $G$. If one
expresses this set as a union of $G$-orbits $\bigcup_PGP$, each counted with
the multiplicity $m_P$ of $z_P$ as a root of $f(z)=0,$ then the polynomial
$\prod_Pf_P^{m_P}(z)$ has the same zeros and multiplicity as $f(z)$, hence
differs from $f$ by a constant factor. \ As form on $V$, $f=$
$\prod_Pf_P^{m_P}$ is an absolute invariant if the multiplicative character
$\prod_P\mu_P^{m_P}$ is $=1$. In the case when $G$ is a finite subgroup
$K=\langle p_A,p_B,p_C\rangle$ of \text{SU}$(2)$ one obtains generators for
$R$ in this way by restricting $P$ to run over the eigenspace $P_A,P_B,P_C$
of the three maximal abelian subgroups $T_A,T_B,T_C$of $G$. In case $G$ is
the reflection group $K'=\langle p_A,p_B,p_C\rangle'$ two of such invariants
suffice as algebraically independent generators for $R$. The decomposition
$S\mathfrak{f}=\prod\mathfrak{F}_g^e$ of an ideal $\mathfrak{f}=Rf$ generated
by such an invariant $f=$ $\prod_Pf_P^{m_P}$ is clear. 

A question arises at this point. The discussion here, and the proof of the
theorem by reduction to ramification of primes in a Dedekind domain, depends
heavily on $\dim V=2$ . The question is whether the decomposition
$S\mathfrak{p}=\prod\mathfrak{P}_g^e$ nevertheless remains valid if $G$ is a
reflection group in arbitrary dimension and $\mathfrak{P}$ the ideal
vanishing on an arbitrary intersection of reflecting hyperplanes.

\medskip

\noindent\textbf{Geometry}. The discussion of the ramification of primes in
the ring of invariants proceeded in the language of commutative algebra.
Another perspective emerges if it is rephrased in the language of algebraic
geometry. The translation runs as follows. The rings $S$ and $R$ are the
regular functions on the affine varieties $V$ and $Q:=G\backslash V$, the
inclusion $S\supset R$ being the pull-back via the quotient map
$\pi:V\rightarrow Q$. These varieties come with actions of the multiplicative
group $\mathbb{C}^\times$ whose quotients give projective varieties
$\bar{V}\approx\mathbb{P}^1$ and $\bar{Q}$ together with a map
$\bar{\pi}:\bar{V}\rightarrow\bar{Q}$. $V$ and $Q$ are affine algebraic
surfaces, $\bar{V}$ and $\bar{Q}$ projective curves. The maps $\pi$ and
$\bar{\pi}$ are ramified coverings.

A 1-dimensional subspace $P$ of $V$ corresponds to a point $\bar{P}$ of
$\bar{V}$ and its image $p$ in $Q$ to a point $\bar{p}$ in $\bar{Q}$. If
$\mathfrak{p}$ is the ideal of $R$ defining the image $p=\pi(P)$ of $P$ in
$Q$, then $S\mathfrak{p}$ is the ideal of $S$ defining the inverse image
$\pi^{-1}(p)=GP=\bigcup gP$ in $V$. The factorization
$S\mathfrak{p}=\prod_{g\in I\backslash G}\mathfrak{P}_g^e$ corresponds to the
decomposition of $\pi^{-1}(p)$ into the individual lines $gP$ each counted
with multiplicity $e$. $S/\mathfrak{p}$ is the ring of 'functions' on
$\pi^{-1}(p)$ and $S/\mathfrak{p}=\prod_{g\in I\backslash
G}S/\mathfrak{P}_g^e$ is its decomposition into functions supported on the
individual lines $gP$. If $P$ is the replaced by a point $x$ in it, then
$S\mathfrak{m}=\prod_{g\in I\backslash G}\mathfrak{M}_g^e$ defines the
decomposition of the fiber $\pi^{-1}\pi(x)=Gx$ into multiple points and
$S/\mathfrak{m}=\prod_{g\in I\backslash G}S/\mathfrak{M}_g^e$ the
decomposition of the ring of 'functions' on the orbit $Gx$into functions
supported at its individual points $gx$. \ The isotropy group $I$ of the
base-point $x$ is normal in the group $H$ leaving $P$ invariant and the
groups $G$ and $H$ act on opposite sides on the orbit: if $G$ acts on
$Gx\approx G/I$ on the left, then $H$ acts on right via $(gx)h:=ghx$. The
difficulty with this picture comes from the fact that the points $gx$ in an
orbit may have higher multiplicity, accounting for the fact that the
'functions' take on values in the rings $S/\mathfrak{M}_g^e$ depending on the
points, rather than in the field $S/\mathfrak{M}_g=\mathbb{C}$. It is for
this reason that the algebraic description seems preferable.

Nevertheless, the interpretation of ramification as representing fiber
decomposition with multiplicity, i.e.
\begin{equation*}
S\mathfrak{m}=\prod\mathfrak{M}_g^e\quad{\text{represents}}\quad\pi^{-1}[%
\pi(x)]=\sum e[x_g],
\end{equation*}
is highly suggestive: it suggests a realization of the representation of $G$
with character $\pi_i$ by monodromy in this fiber system, with $G$ acting on
a fiber of the map $\pi:V\rightarrow Q$ by continuous transport around loops
in the set of regular values. Of particular interest in this connection is
the \emph{different} ideal of $S$ over $R$ \lbrack Zariski and Samuel, 1958,
Ch.V, \S 11\rbrack. In geometric terms it is the ideal $\mathfrak{D}$ of $S$
defining the ramification locus, i.e. the inverse image of the set of
critical values of $\pi:V\rightarrow Q$. Equivalently, $\mathfrak{D}$ is the
ideal in the homogenous coordinate ring $S$ of $\bar{V}\approx\mathbb{P}^1$
defining the ramification locus of the ramified covering
$\bar{\pi}:\bar{V}\rightarrow\bar{Q}$ of compact Riemann surfaces. It factors
as $\mathfrak{D}=\prod_P\mathfrak{P}^{m_P}$ where $\bar{P}$ runs over the
points of $\bar{V}$ at which the Jacobian $D$ of $\bar{\pi}$ vanishes with
multiplicity $m_P\neq 0$. It is in fact the principal ideal
$\mathfrak{D}=SD$. In case $G$ is a reflection group (in arbitrary dimension)
the ramification locus is the arrangement of reflecting hyperplanes and $D$
the product of the linear forms defining them. \lbrack Bourbaki, 1968,
p.116\rbrack. 

In case $G$ is a finite subgroup of \text{U}$(2)$, e.g. a Klein group, the
ramification locus is the union of the singular lines $P$ which have the
maximal abelian subgroups $H$ of $G$ as stabilizers. As mentioned in \S 2,
each $H$ is the stabilizer of two singular lines $P,P^{-}$ which may or may
not be interchanged by the group $N(H)/H$, of order 2 or 1.
$D\propto\prod_P(\alpha_P)^{m_P}$ is a constant multiple of the product of
linear forms $a_P$ defining the singular lines, $a_P$ being taken with
multiplicity $m_P=|\bar{H}\mathopen|-1$, $\bar{H}$ $:=$ the image of $H$ in
\text{$\bar{\text{U}}$}$(V):={\text{U}}(V)/{\text{U}}(\mathbb{C})$. If the
$\alpha_P$s are grouped into $G$-orbits the product becomes 
\begin{equation*}
D\propto f_A^{p_A-1}f_B^{p_B-1}\cdots
\end{equation*}
Here $H_A,H_B,\cdots$ runs over a complete system of representatives for the
conjugacy classes of $H$s, $f_A=\prod_{g\in[G/H_A]}(a_Ag),$ and
$p_A=|\bar{H}_A|$. 

On the other hand, the Jacobian is
$D\propto\partial(f_1,f_2)/\partial(z_1,z_2)$ for any are linear coordinates
$z_1,z_2$ of $V$ and any two generic forms $f_1,f_2$ of the type $f_P$, for
which $\bar{H}=\{1\}$. If one compares degrees in the relation
$\partial(f_1,f_2)/\partial(z_1,z_2)\propto f_A^{p_A-1}f_B^{p_B-1}\cdots$ one
finds
\begin{equation*}
2(|\bar{G}\mathopen|-1)=\frac{2|\bar{G}|}{|N(\bar{H})/|\bar{H}|}(|\bar{H}%
\mathopen|-1)+\cdots.
\end{equation*}
One thus returns to the class equation for subgroups of \text{U}$(2)$, the
starting point for their classification. The discussion is based on \lbrack
Klein, 1884, Ch.5, \S 2\rbrack.


\bigskip

\begin{thebibliography}{15}


\bibitem{}N. Bourbaki, Alg\a`ebre commutative. Chapitres 1 et 2. Hermann,
Paris, 1961.


\bibitem{}N. Bourbaki, Groupes et alg\a`ebres de Lie. Chapitres 4,5 et 6.
Hermann, Paris, 1968.

\bibitem{}N. Bourbaki, Alg\a`ebre commutative. Chapitres 5,6 et 7. Hermann,
Paris, 1985.


\bibitem{}H.S.M. Coxeter, Regular Complex Polytopes. Second Edition.
Cambridge University Press, 1974.


\bibitem{}G. Gonzales-Sprinberg and J.-L. Verdier, Construction
g\a'eom\a'etrique de la corrrespondence de McKay. Ann. Sci. Ecole Norm. Sup.
t.16, \text{n}$^{\text{o}}3$, 410-449(1983).


\bibitem{}E. Hecke. Vorlesungen \"uber die Theorie der Algebraischen Zahlen.
Leipzig 1923. Reprinted by Chelsea Publishing Co., New York, 1970.


\bibitem{}F. Klein. Vorlesungen \"uber das Ikosaeder und die Aufl\"osung der
Gleichung vom f\"unften Grade. Teubner, Leipzig, 1884.


\bibitem{}H. Kn\"orrer. Group representations and the resolution of rational
double points. In: Finite groups - Coming of Age. Proceedings, Montreal 1982
(J. Mckay, ed.). Contemporary Math., v.45, AMS, 175-222, Providence, 1985.


\bibitem{}B. Kostant, The McKay correspondence, the Coexeter element, and
representation theory. In: \a'Elie Cartan el les math\a'emaiques
d'aujourd'hui (Lyon, 1984). Ast\a'erisque, Hors s\a'eries, 209-255 (1985).


\bibitem{}G. Lusztig, Some examples of square integrable representations of
semisimple p-adic groups. Trans. AMS 277, 153-215 (1983).


\bibitem{}G. Lusztig, Subregular nilpotent elements and bases in K-theory.
Cand. J. Math. 51(6), 1194-1225 (1999).


\bibitem{}J. McKay. Graphs, singularities, and finite groups. AMS, Proc.
Symp. Pure Math. Vol. 37, 183-186, 1980.


\bibitem{}J.P. Serre. Linear Represenations of Finite Groups. Springer
Verlag, New York, 1977.


\bibitem{}G.C. Shephard and J.A. Todd. Finite unitary reflection groups. Can.
J. Math. 6, 111-135 (1954).


\bibitem{}T.A. Springer, Poincar\a'e series of binary polyhedral groups and
McKay's correspondence. Math. Ann. 278 , 587-598 (1985).


\bibitem{}R. Steinberg, Finite subgroups of \text{SU}$_2$, affine Dynkin
diagrams and affine Coxeter elements. Pac. J. Math. 118, 587-598 (1985).
Preprint 1982.


\bibitem{}H. Weyl. The Classical Groups. Their invariants and
Representations. Second Edition. Princeton University Press, Princeton, 1939.


\bibitem{}O. Zariski and P. Samuel. Commutative Algebra I. D. Van Nostrand
Company, Inc. Princeton, 1958.

\end{thebibliography}
\end{document}